\documentclass[12pt]{amsart}
\usepackage{amsmath}
\usepackage{amssymb}
\usepackage{latexsym}
\def\ds{\displaystyle}
\newtheorem{theorem}{Theorem}
\newtheorem{conjecture}[theorem]{Conjecture}

\newtheorem{lemma}[theorem]{Lemma}
\newtheorem{proposition}[theorem]{Proposition}
\newtheorem{definition}[theorem]{Definition}

\newcommand{\C}{{\mathbb{C}}}

\newcommand{\Ker}{\operatorname{Ker}}

\title[Discontinuity of the Lempert function of the spectral ball]
{Discontinuity of the Lempert function and the Kobayashi--Royden
metric of the spectral ball}
\author{Nikolai Nikolov, Pascal J. Thomas, W\l odzimierz Zwonek}

\address
{Institute of Mathematics and Informatics\\ Bulgarian Academy of
Sciences\\ Acad. G. Bonchev 8, 1113 Sofia,
Bulgaria}\email{nik@math.bas.bg}

\address{Laboratoire Emile Picard, UMR CNRS 5580\\
Universit\'e Paul Sabatier, 118 Route de Narbonne\\ F-31062 Toulouse
Cedex, France} \email{pthomas@cict.fr}

\address{Instytut Matematyki, Uniwersytet Jagiello\'nski, Reymonta 4,
30-059 Krak\'ow, Poland}\email{Wlodzimierz.Zwonek@im.uj.edu.pl}



\begin{document}

\begin{thanks}{This work was initialized during the stay of the first and
second named authors at the Jagiellonian University, Krak\'ow in
October, 2006, supported by the EGIDE program. They wish thank
both institutions. The third author was supported by the KBN research grant No. 1 PO3A 005 28.}
\end{thanks}

\begin{abstract} Some results on the discontinuity properties of the Lempert function and the
Kobayashi pseudometric in the spectral ball are given.
\end{abstract}

\maketitle

\section{Introduction and results}

Let $\mathcal M_n$ be the set of all $n\times n$ complex matrices.
For $A\in\mathcal M_n$ denote by $sp(A)$ and $\ds
r(A)=\max_{\lambda\in sp(A)}|\lambda|$ the spectra and the
spectral radius of $A,$ respectively. The spectral ball $\Omega_n$
is the set
$$\Omega_n=\{A\in\mathcal M_n:r(A)<1\}.$$
The Nevanlinna--Pick problem in $\Omega_n$ (or the spectral
Nevanlinna--Pick problem) is the following one: given $N$ points
$a_1,\dots,a_N$ in the unit disk $\Bbb D\subset\Bbb C$ and $N$
matrices $A_1,\dots, A_N\in\Omega_n$ decide whether there is a
holomorphic map $\varphi \in\mathcal O(\Bbb D,\Omega_n)$ such that
$\varphi(a_j )=A_j$, $1\le j \le N$.  This problem has been
studied by many authors; we refer the reader to
\cite{Agl-You1,Agl-You2,Agl-You3,BFT,Cos2,Cos3} and the references
there.

The study of the spectral Nevanlinna--Pick problem in the case $N=2$
reduces to the computation of the Lempert function of the spectral
ball. Recall that for a domain $D\subset\Bbb C^m$ the Lempert function of the domain $D$ is defined as
follows:
$$ l_D(z,w):=\inf\{|\alpha|:\exists\varphi\in\mathcal O(\Bbb
D,D):\varphi(0)=z,\varphi(\alpha)=w\},\; z,w\in D.$$

The infinitesimal version of the above problem, the so-called
spectral Carath\'eodory--Fej\'er problem is the following one:
given $N+1$ matrices $A_0,\dots, A_N$ in $\mathcal M_n$ decide
whether there is a map $\varphi\in \mathcal O(\Bbb D,\Omega_n)$
such that $A_j=\varphi^{(j)}(0),$ $0\le j\le N$. This problem has
been studied in \cite{HMY}.

The study of the spectral Carath\'eodory-Fejer problem in the case
$N=1$ reduces to the computation of the Kobayashi--Royden
pseudometric of the spectral ball. Recall that for a domain
$D\subset\Bbb C^m$ the Kobayashi--Royden pseudometric is defined
as follows:
\begin{align*}
\kappa_D(z;X):=\inf\{|\alpha|:\exists\varphi\in&\mathcal O(\Bbb
D,D):\\&\varphi(0)=z,\alpha\varphi'(0)=X\},\; z\in D,X\in\Bbb C^m.
\end{align*}

In this note we point out some of the non-stability phenomena of both spectral
problems
which complicate their study.

First note that if we replace each of the matrices in the spectral
Nevanlinna--Pick problem by similar ones then we do not change its
solution.\footnote{Indeed, assume that $\varphi\in\mathcal O(\Bbb
D,\Omega_n),$ $\varphi(a_j )=A_j$ and $A_j\sim\tilde A_j,$ $1\le j
\le N$. Then $\tilde A_j=e^{\hat A_j}A_j e^{-\hat A_j}$ for some
$\hat A_j\in\mathcal M_n.$ In the standard way we may find
$\hat\varphi\in\mathcal O(\Bbb D,\mathcal M_n)$ with
$\hat\varphi(a_j)=\hat A_j,$ $1\le j\le N.$ Then $\tilde\varphi:=
e^{\hat\varphi}\cdot \varphi\cdot e^{-\hat\varphi}\in\mathcal
O(\Bbb D,\Omega_n)$ and $\tilde\varphi(a_j )=A_j.$} A natural
reduction of the problems is then to associate to each matrix its
spectrum, or, in order to deal with $n$-tuples of complex numbers,
the coefficients of its characteristic polynomial
$$
P_A(t):= \det(tI-A)=t^n+\sum_{j=1}^n(-1)^j\sigma_j (A)t^{n-j},
$$
where $I\in \mathcal M_n$ is the unit matrix,
$$\sigma_j(A):=\sigma_j(\lambda_1,\dots,\lambda_n):=\sum_{1\le k_1<\dots<k_j\le
n}\lambda_{k_1}\dots,\lambda_{k_j}$$ and $\lambda_1,\dots,\lambda_n$
are the eigenvalues of $A.$

Put $\sigma=(\sigma_1,\dots,\sigma_n).$ We shall consider $\sigma$
as a map either from  $\mathcal O(\mathcal M_n,\Bbb C^n),$ or from
$\mathcal O(\Bbb C^n,\Bbb C^n).$ The set $$\Bbb
G_n:=\{\sigma(A):A\in\Omega_n\}$$ is called the symmetrized
$n$-disk, and has been widely studied; we refer the reader to
\cite{Agl-You3,Cos3,Jar-Pfl1,NPTZ,NPZ} and references there.

We recall a few definitions from linear algebra.
\begin{definition}
Given a matrix $A\in\mathcal M_n$, the \emph{commutant} of $A$ is
$$
\mathcal C (A) := \{ M \in\mathcal M_n : MA=AM\},
$$
and the set of polynomials in $A$, $ \mathcal P (A) \subset \mathcal C (A)$ is
given by
$$
\mathcal P (A) := \{ M \in\mathcal M_n : M=p(A), \mbox{ for some } p\in \mathbb C[X]\}.
$$
\end{definition}
\begin{definition}
Given $(a_0, \ldots, a_{n-1}) \in \C^n$, the associated \emph{companion matrix} is
$$
\left(
\begin{array}{ccccc}
0&&&&-a_0\\
1&0&&&\vdots\\
&1&\ddots&&\vdots\\
&&\ddots&0&-a_{n-2}\\
&&&1&-a_{n-1}
\end{array}
\right) .
$$
The companion matrix associated to a matrix $A$ is the one associated
to the coefficients of its characteristic polynomial, namely we set $a_j= (-1)^{n-j}\sigma_{n-j}(A)$,
$0\le j \le n-1$.
\end{definition}

\begin{proposition}
\label{equiv}
A matrix $A\in\mathcal M_n$ with the following equivalent properties
is called \emph{non-derogatory}.
\begin{enumerate}

\item $A$ is similar to its companion matrix.

\item There exists a cyclic vector for $A.$

\item The characteristic and minimal polynomials of $A$ coincide.

\item Different blocks in the Jordan normal form of $A$ correspond
to different eigenvalues (that is, each eigenspace is of dimension
exactly $1$).

\item $\mathcal C(A) = \mathcal P(A)$.

\item $rank(\sigma_{\ast,A})=n.$

\item $\dim \mathcal C(A) = n$.

\item If $\Phi_A: \mathcal M_n^{-1 } \longrightarrow \mathcal
M_n$, where $ \mathcal M_n^{-1 }$ stands for the set of invertible
matrices, is defined by $\Phi_A(P):=P^{-1} A P$, then
$rank((\Phi_A)_{\ast,I_n})\\=n^2-n$ (its maximal possible value).
\end{enumerate}
\end{proposition}

Most of those properties can be found in \cite[pp.
135--147]{HoJo}; more precise references and (easy) complements
are given in Section \ref{appendix}.

Recall that the $j$-th coordinate of $\sigma_{\ast,A}(B)$ is the sum
of all $j\times j$ determinants obtained by taking a principal
$j\times j$ submatrix of $A$ and replacing one column by the
corresponding entries of $B.$ In particular, the first coordinate
of $\sigma_{\ast,A}(B)$ equals $tr(B).$

Denote by $\mathcal C_n$ the set of all non-derogatory matrices in
$\Omega_n.$ Obviously $\mathcal C_n$ is an open and dense subset
of $\Omega_n$.

Note that if $A_1,\dots,A_N$ belong to $\mathcal C_n$, then
any mapping $\varphi\in\mathcal O(\Bbb b D,\Bbb G_n)$ with
$\varphi (\alpha_j)=\sigma(A_j)$ can be lifted to a mapping
$\tilde\varphi\in\mathcal O(\Bbb D,\Omega_n)$ with $\tilde\varphi
(\alpha_j)=A_j,$ $1\le j\le N$ (see \cite{Agl-You1}). This means
that in a generic case the spectral Nevanlinna--Pick problem for
$\Omega_n$ with dimension $n^2$ can be reduced to the standard
Nevanlinna--Pick problem for $\Bbb G_n$ with dimension $n.$

As a consequence of the existence of the lifting above, we have
the equality
\begin{equation}\label{1}
l_{\Omega_n}(A,B)=l_{\Bbb G_n}(\sigma(A),\sigma(B)),\quad
A,B\in\mathcal C_n.
\end{equation}
Note that $\Bbb G_n$ is a taut domain (cf. \cite{Edi-Zwo},
\cite{Jar-Pfl1}). In particular, there always exist extremal discs
for $l_{\Bbb G_n}$ and $l_{\Bbb G_n}$ is a continuous function. Thus
the spectral Nevanlinna--Pick problem with data
$(\alpha_1,A_1),(\alpha_1,A_2)\in\Bbb D\times\mathcal C_n$ is
solvable if and only if
$$l_{\Bbb G_n}(\sigma(A),\sigma(B))\le m(\alpha_1,\alpha_2)
:=\left|\frac{\alpha_1-\alpha_2}{1-\alpha_1\overline{\alpha_2}}\right|.$$

An explicit formula for $l_{\Bbb G_2}$ is found in
\cite{Agl-You3}. The proof there is based on studying the complex
geodesics of $\Bbb G_2.$ It turns out that $\tanh^{-1}l_{\Bbb
G_2}$ coincides with the Carath\'eodory distance of $\Bbb G_2$. On
the other hand, $\Bbb G_2$ cannot be exhausted by domains
biholomorphic to convex domains (see \cite{Cos1}, \cite{Edi}). So
$\Bbb G_2$ serves as the first counterexample to converse of the
Lempert theorem (cf. \cite{Jar-Pfl1}). In spite of this
phenomenon, $\tanh^{-1}l_{\Bbb G_n},$ does not even satisfy the
triangle inequality for $n>2$, that is, it does not coincide with
the Kobayashi distance of $\Bbb G_n$ for $n>2$ (see \cite{NPZ}).

The behavior of $l_{\Omega_n}$ is much more complicated when one of
the arguments is derogatory. However, if $A$ is a scalar matrix, say
$A=tI,$ $t\in\Bbb D,$ then (cf. \cite{Agl-You1})
\begin{equation}\label{e2}
l_{\Omega_n}(tI,B)=\max_{\lambda\in sp(B)}m(t,\lambda).
\end{equation}
To prove (\ref{e2}), observe first that
$B\to(B-tI)(I-\overline{t}B)^{-1}$ is an automorphism of
$\Omega_n.$ So we may assume that $t=0.$ Then it remains to make
use of the fact that $l_{\Omega_n}(0,B)$ equals the Minkowski
function of the balanced domain $\Omega_n$ at $B,$ that is,
$r(B).$

Since the $2\times2$ derogatory matrices are scalar, we also get
that the function $l_{\Omega_2}$ is not continuous at the point
$(A,B)$ if and only if one of the matrices is scalar, say $A,$ and
the other one has two distinct eigenvalues (see \cite{Cos2}). In
this case even $l_{\Omega_2}(\cdot,B)$ is not continuous at $A$
(but $l_{\Omega_2}(A,\cdot)$ is continuous at $B$). We shall show
that this phenomenon extends to $\Omega_n.$

\begin{proposition}\label{2} For $B\in\mathcal C_n$ and $t\in\Bbb
D$ the following conditions are equivalent:

(i) the eigenvalues of $B$ are equal;

(ii) the function  $l_{\Omega_n}$ is continuous at the point
$(tI,B);$

(iii) the function $l_{\Omega_n}(\cdot,B)$ is continuous at the
point $tI;$
\end{proposition}

At the infinitesimal level of the Kobayashi--Royden pseudometric,
for $A\in\mathcal C_n$ and $B\in\mathcal M_n$ one has that (see
Theorem 2.1 in \cite{HMY})
\begin{equation}\label{e3}
\kappa_{\Omega_n}(A;B)=\kappa_{\Bbb
G_n}(\sigma(A);\sigma_{*,A}(B)).
\end{equation}
Since $\kappa_{\Omega_n}(A;B)\ge \kappa_{\Bbb
G_n}(\sigma(A);\sigma_{*,A}(B))$ for
$(A;B)\in\Omega_n\times\mathcal M_n,$ $\kappa_{\Omega_n}$ is an
upper semicontinuous function and $\kappa_{\Bbb G_n}$ is
continuous (because $\Bbb G_n$ is a taut domain), we get that
$\kappa_{\Omega_n}$ is a continuous function at any point
$(A;B)\in\mathcal C_n\times \mathcal M_n.$

The things are more complicated if $A\not\in\mathcal C_n.$

\begin{proposition}\label{3} For $B\in\mathcal M_n$ and $t\in\Bbb D$
the following conditions are equivalent:

(i) the eigenvalues of $B$ are equal;

(ii) the function  $\kappa_{\Omega_n}$ is continuous at the point
$(tI;B);$

(iii) the function $\kappa_{\Omega_n}(\cdot;B)$ is continuous at
the point $tI;$
\end{proposition}

Note that, similarly to the equality (\ref{e2}), one has that
\begin{equation}\label{e4}
\kappa_{\Omega_n}(tI;B)=\frac{\max_{\lambda\in
sp(B)}|\lambda|}{1-|t|^2}.
\end{equation}

As a consequence of our considerations, we may also identify in a
simple way the convex hull $\hat\Omega_n$ of $\Omega_n$.

\begin{proposition}\label{4}
$\hat\Omega_n=\{ A \in\mathcal M_n:|tr(A)|<n\}.$
\end{proposition}

We now turn to analyzing the failure of hyperbolicity of
$\Omega_n.$ Observe first that if $sp(A)\neq sp(B)$ then
$\sigma(A)\neq\sigma(B)$ and hence
$$l_{\Omega_n}(A,B)\ge l_{\Bbb G_n}(\sigma(A),\sigma(B))>0.$$
Then as a consequence of the proof of Lemma 13 in \cite{Edi-Zwo} we
have the following

\begin{proposition}\label{5} For any $A,B\in\Omega_n$ the equality
$l_{\Omega_n}(A,B)=0$ holds if and only if $sp(A)=sp(B)$. Moreover, in
this case there is a $\varphi\in\mathcal O(\Bbb C,\Omega_n)$ with
$\varphi(0)=A,$ $\varphi(1)=B$ and $sp(\varphi(\lambda))=sp(A)$
for any $\lambda\in\Bbb C.$
\end{proposition}

It is natural to consider the infinitesimal version of this
proposition.

First, note that the equality (\ref{e4}) implies that if
$A\in\Omega_n$ is a scalar matrix and $B\in\mathcal M_n,$ then
$\kappa_{\Omega_n}(A;B)=0$ if and only if $sp(B)=0.$ Conversely,
$A$ is scalar and $sp(B)=0,$ then the linear mapping $p:\lambda\to
A+\lambda B$ has the following properties: $p(0)=A,$ $p'(0)=B$ and
$sp(p(\lambda))=sp(A)$ for any $\lambda\in\Bbb C.$

On the other hand, the equality (\ref{e3}) implies that if
$A\in\mathcal C_n$ and $B\in\mathcal M_n,$ then
$\kappa_{\Omega_n}(A;B)=0$ if and only if $\sigma_{\ast,A}(B)=0.$
Moreover, the following is true.

\begin{proposition}
\label{nonderent} If $A\in\mathcal C_n,$ $B\in\mathcal M_n$ and
$\kappa_{\Omega_n}(A;B)=0,$ then there is a mapping $\varphi \in
\mathcal O(\C,\Omega_n)$ with $\varphi(0)=A,$ $\varphi'(0)=B$ and
$sp(\varphi(\lambda))=sp(A)$ for any $\lambda\in\C.$
\end{proposition}

These observations let us state the following

\begin{conjecture}\label{6} If $A\in\Omega_n,$ $B\in\mathcal
M_n$ and $\kappa_{\Omega_n}(A;B)=0,$ then there is a polynomial
mapping $p:\Bbb C\to\Omega_n$ of degree at most $n$ with $p(0)=A,$
$p'(0)=B$ and $sp(p(\lambda))=sp(A)$ for any $\lambda\in\Bbb C.$
\end{conjecture}

To support this conjecture, we shall prove it for $n=2.$
\smallskip

The rest of the paper is organized as follows. The proof of
Proposition \ref{2} is given in Section \ref{2}. Section \ref{3}
contains the proofs of Propositions \ref{3}, \ref{4} and
\ref{nonderent}, as well as Conjecture \ref{6} for $n=2.$ The
proof of Proposition \ref{equiv} is discussed in Section 4.

\section{Proof of Proposition \ref{2}}

We shall need the following

\begin{proposition}\label{7}

(i) If $A,B\in\Omega_n,$ then $$l_{\Bbb G_n}(\sigma(A),\sigma(B))\le
l_{\Omega_n}(A,B)\le\min_{\pi}\max_{1\le j\le n}
m(\lambda_j,\mu_{\pi(j)}),$$ where
$sp(A)=\{\lambda_1,\dots,\lambda_n\},$
$sp(B)=\{\mu_1,\dots,\mu_n\},$ and the minimum is taken over all
permutations $\pi$ of $\{1,\dots,n\}.$

(ii) (see Theorem 5.2 in \cite{Cos3}) If the eigenvalues of
$B(z)\in\mathcal O(\Bbb D,\Omega_n)$ have the form
$e^{i\theta}\frac{1-\overline{\alpha_j}}{1-\alpha_j}
\frac{z-\alpha_j}{1-z\overline{\alpha_j}},$ then
$$l_{\Bbb G_n}(\sigma(B(z)),\sigma(B(w)))=l_{\Omega_n}(B(z),B(w))=m(z,w).$$

(iii) If $B\in\Omega_n$ and $t\in\Bbb D,$ then the eigenvalues of
$B$ are equal if and only if $$l_{\Bbb
G_n}(\sigma(tI),\sigma(B))=\max_{\lambda\in sp(B)}m(t,\lambda).$$
\end{proposition}

\noindent{\bf Remark.} One may conjecture that Proposition \ref{7}
(ii) describes all the possibilities for the equality
$$l_{\Bbb G_n}(\sigma(A),\sigma(B))=\min_{\lambda\in sp(A)}\max_{\mu\in
sp(B)}m(\lambda,\mu).$$

Assuming Proposition \ref{7} (iii), we are ready to prove
Proposition \ref{2}.

The implication $(ii)\Rightarrow(iii)$ is trivial. For the rest of
the proof we may assume that $t=0.$

We shall show that $(i)\Rightarrow(ii)$ for any $B\in\Omega_n.$
Let $(A_j)\to 0$ and $(B_j)\to B.$ Then, by Proposition \ref{7}
(iii) and (\ref{e2}),
$$l_{\Omega_n}(A_j,B_j)\ge l_{\Bbb
G_n}(\sigma(A_j),\sigma(B_j))\to l_{\Bbb
G_n}(0,\sigma(B))=r(B)=l_{\Omega_n}(0,B).$$ Thus the function
$l_{\Omega_n}(\cdot,B)$ is lower semicontinuous at the point
$(0,B).$ Since it is (always) upper semicontinuous, we conclude
that it is continu\-ous at this point.

It remains to prove that  $(iii)\Rightarrow(i).$ Since $\mathcal
C_n$ is a dense subset in $\Omega_n,$ we may find $\mathcal
C_n\supset(A_j)\to 0.$ Then, by (\ref{e2}) and (\ref{1}),
$$r(B)=l_{\Omega_n}(0,B)\leftarrow
l_{\Omega_n}(A_j,B)=l_{\Bbb G_n}(\sigma(A_j),\sigma(B))\to l_{\Bbb
G_n}(0,\sigma(B))$$ and hence $l_{\Bbb G_n}(0,\sigma(B))=r(B).$ It
follows by Proposition \ref{7} (iii) that the eigenvalues of $B$
are equal.

This completes the proof of Proposition \ref{2}.
\smallskip

\noindent{\bf Proof of Proposition \ref{7} (i).} The first
inequality is trivial.

To prove the second one recall that (see the footnote on page 2)
$$l_{\Omega_n}(A,B)=l_{\Omega_n}(A',B'),\quad A'\sim A,\ B'\sim B.$$
So we may assume $A=(a_{jk})$ and $B=(b_{jk})$ are Jordan matrices
with $$\max_{1\le j\le n}m(a_{jj},b_{jj})=s:= \min_{\pi}\max_{1\le
j\le n} m(\lambda_j,\mu_{\pi(j)}).$$ Let $s_1> s.$ Then we may
choose $\varphi_{jj}\in\mathcal O(\overline{\Bbb D},\Bbb D)$ such
that $\varphi_{jj}(0)=a_{jj}$ and $\varphi_{jj}(s_1)=b_{jj}.$ For
$\zeta\in\Bbb C$ set
$$\varphi_{jk}(\zeta)=\left\{\begin{array}{ll}
0,&j>k\\
a_{jk}+\frac{b_{jk}-a_{jk}}{s_1}\zeta,&j<k.
\end{array}\right.$$
Now $\varphi=(\varphi_{jk})\in\mathcal O(\overline{\Bbb
D},\Omega_n)$ which shows that $l_{\Omega_n}(A,B)<s_1.$ Since
$s_1>s$ was arbitrary, we are done.
\smallskip

\noindent{\bf Remark.} Obvious modifications in the above proof
imply Proposition \ref{5}.
\smallskip

\noindent{\bf Proof of Proposition \ref{7} (iii).} If the
eigenvalues of $B$ are equal, say to $\lambda,$ then
$$l_{\Bbb G_n}(\sigma(tI),\sigma(B))=m(t,\lambda)$$
by Proposition \ref{7} (ii) with $\alpha_j=0$ (or, directly,
considering the mapping $\zeta\to\sigma(\zeta,\dots,\zeta)$ shows
that second inequality in Proposition \ref{7} (i) becomes
equality).

To prove the converse, we shall need the following

\begin{lemma}\label{8} Let $\varepsilon_1,\dots,\varepsilon_n\in\Bbb T
=\partial\Bbb D$ be pairwise different points. Then for any
$\lambda_1,\dots,\lambda_n\in\Bbb D,$ there are $\beta\in\Bbb D$
and a Blaschke product $\mathcal B$ of order $\le n$ with
$\mathcal B(0)=0,\mathcal
B(\varepsilon_1\beta)=\lambda_1,\dots,\mathcal
B(\varepsilon_n\beta)=\lambda_n.$
\end{lemma}

Assuming Lemma \ref{8}, we shall complete the proof of Proposition
\ref{7} (iii). We may assume that $t=0.$ Let
$\lambda_1,\dots,\lambda_n$ be the eigenvalues of $B.$ Set $\root
n\of 1=\{\varepsilon_1,\dots,\varepsilon_n\}.$ Let $\beta\in\Bbb
D$ and $\mathcal B$ be as in Lemma \ref{8}. Consider the mapping
$$\zeta\to f_{\mathcal B}(\zeta):=\sigma(\mathcal B(\varepsilon_1
\root n\of {\zeta}),\dots,\mathcal B(\varepsilon_n\root
n\of{\zeta}))$$ (where $\root n\of{\zeta}$ is arbitrary chosen). It
is easy to see that $f_{\mathcal B}\in\mathcal O(\Bbb D,\Bbb G_n).$
Hence $l_{\Bbb G_n}(0,\sigma(B))\le|\beta|^n.$ It remains to prove
that if $\ds|\beta|^n\ge\max_{1\le j\le n}|\lambda_j|$ then
$\lambda_1=\dots=\lambda_n.$ We may assume that
$$\mathcal B(\zeta)=z\frac{a_0\zeta^k+a_1\zeta^{k-1}+\dots+a_k}{\bar a_k\zeta^k+\bar
a_{k-1}\zeta^{k-1}+\dots+\bar a_0},$$ where $a_0=1$ and $k\le
n-1.$ Then for $1\le j\le n $ one has that
$$|\beta|^{n-1}|\overline{a_k}(\varepsilon_j\beta)^k+\overline{
a_{k-1}}(\varepsilon_j\beta)^{k-1}+\dots+\overline{a_0}|\ge
|a_0(\varepsilon_j\beta)^k+a_1(\varepsilon_j\beta)^{k-1}+\dots+a_k|.$$
Squaring both sides of this inequality, we get that
$$|\beta|^{2n-2}(\sum_{s=0}^k|a_s|^2|\beta|^{2s}+2\Re\sum_{0\le
p<s\le
k}a_p\overline{a_s}\beta^s\overline{\beta}^p\varepsilon_j^{s-p})\ge$$
$$\sum_{s=0}^k|a_s|^2|\beta|^{2(k-s)}+2\Re\sum_{0\le p<s\le
k}a_p\overline{a_s}\beta^{k-p}\overline{\beta}^{k-s}\varepsilon_j^{s-p}.$$
Summing these inequalities for $j=1,\dots,n$ we get that
$$|\beta|^{2n-2}\sum_{s=0}^k|a_s|^2|\beta|^{2s}\ge
\sum_{s=0}^k|a_s|^2|\beta|^{2(k-s)},$$ that is,
$$\sum_{s=0}^k|a_s|^2(|\beta|^{2(n+s-1)}-|\beta|^{2(k-s)}|)\ge
0.$$ Since $k\le n-1,$ then $k-s<n+s-1$ if $s>0$ and hence
$a_s=0.$ On the other hand, $a_0\neq 0$ and thus $k=n-1.$ It
follows that $\mathcal B(z)=z^n$ and then
$\lambda_1=\dots=\lambda_n.$
\smallskip

\noindent {\bf Proof of Lemma \ref{8}.} Let $S$ be the set of all
$\beta\in\Bbb D$ for which the Nevanlinna--Pick problem with data
$(0,0),(\varepsilon_1\beta,\lambda_1),\dots,(\varepsilon_n\beta,\lambda_n)$
is solvable. We may assume that $0\not\in S;$ otherwise, the
identity mapping does the job. Then we must have $f\in\mathcal
O(\Bbb D,\Bbb D)$ with $\ds
f(\varepsilon_j\beta)=\lambda_j':=\frac{\lambda_j}{\varepsilon_j\beta},$
$1\le j\le n.$ The existence of such a function is equivalent to
the semi-positivity of the matrix
$A(\beta)=[a_{j,k}(\beta)]_{j,k=1}^n,$ where $\ds
a_{j,k}(\beta)=\frac{1-\lambda_j'\bar\lambda_k'}{1-\varepsilon_j
\overline{\varepsilon_k}|\beta|^2}.$ Observe that
$a_{j,k}(\cdot),$ $j\neq k,$ is bounded on $\Bbb D.$ On the other
hand, $\ds\lim_{\beta\to\Bbb T}a_{j,j}(\beta)=+\infty.$ Thus the
matrix $A(\beta)$ is (strictly) positive for $\beta$ near $\Bbb
T.$ Since $0\not\in S,$ it follows that $S$ is a proper non-empty
(circular) closed subset of $\Bbb D.$ So there is a boundary point
$\beta_0\in\Bbb D$ of $S.$ Then $A(\beta_0)$ is not strictly
positive which means that $m=rank(A(\beta_0))$ is not maximal,
that is, $m<n.$ This implies that the respective Nevanlinna--Pick
problem has a unique solution and it is a Blaschke product
$\mathcal{\tilde B}$ of order $m$ (cf. \cite{Gar}). It remains to
set $\mathcal B(\zeta)=\zeta\mathcal{\tilde B}(\zeta).$

\section{Proofs of Propositions \ref{3}, \ref{4}, \ref{nonderent}
and Conjecture \ref{6} for $n=2.$}

\noindent{\bf Proof of Proposition \ref{3}}. The implication
$(ii)\Rightarrow(iii)$ is trivial. For the rest of the proof we
may suppose as above that $t=0.$

We shall show that $(i)\Rightarrow(ii)$ for any $B\in\Omega_n.$
Let $(A_j)\to 0$ and $(B_j)\to B.$ Then, by (\ref{e4}) and the
equality $\ds\kappa_{\Bbb G_n}(0;e_1)=\frac{1}{n}$ (cf.
\cite{NPZ}),
$$\kappa_{\Omega_n}(A_j;B_j)\ge\kappa_{\Bbb
G_n}(\sigma(A_j);\sigma_{\ast,A_j}(B_j))\to \kappa_{\Bbb
G_n}(0;\sigma_{\ast,0}(B))=$$
$$\kappa_{\Bbb G_n}(0;tr(B)e_1)=\frac{|tr(B)|}{n}=r(B)=\kappa_{\Omega_n}(0;B).$$
Thus, the function $\kappa_{\Omega_n}$ is lower semicontinuous at
the point $(0;B).$ Since it is (always) upper semicontinuous, we
conclude that it is continu\-ous at this point.

It remains to prove that  $(iii)\Rightarrow(i).$ Since $\mathcal
C_n$ is a dense subset of $\Omega_n,$ we may find $\mathcal
C_n\supset(A_j)\to 0.$ Then, by (\ref{e4}) and (\ref{e3}),
$$r(B)=\kappa_{\Omega_n}(0;B)\leftarrow
\kappa_{\Omega_n}(A_j;B)=$$ $$\kappa_{\Bbb
G_n}(\sigma(A_j);\sigma_{\ast,A_j}(B))\to\kappa_{\Bbb
G_n}(0;\sigma_{\ast,0}(B))=\frac{|tr(B)|}{n}.$$ Hence $\ds
r(B)=\frac{|tr(B)|}{n},$ that is, the eigenvalues of $B$ are
equal.
\smallskip

\noindent{\bf Proof of Proposition \ref{4}}. Since $\Omega_n$ is a
balanced domain, we have that (cf. \cite{Jar-Pfl1})
$$h_{\hat\Omega_n}=k_{\Omega_n}(0,\cdot),$$
where $h_{\hat\Omega_n}$ and $k_{\Omega_n}$ are the Minkowski
function of $\hat\Omega_n$ and the Kobayashi distance of
$\Omega_n,$ respectively.

On the other hand, since $k_{\Omega_n}$ is a continuous function,
the density of $\mathcal C_n$ in $\Omega_n$ and the equality
(\ref{1}) imply that
$$k_{\Omega_n}(A,B)=k_{\Bbb G_n}(\sigma(A)),\sigma(B)),\quad A,B\in\Omega_n.$$

It follows that for any $t\in\Bbb D\setminus\{0\}$
$$h_{\hat\Omega_n}(A)=\frac{k_{\Omega_n}(0,tA)}{|t|}=$$
$$\frac{k_{\Bbb G_n}(0,\sigma(tA))}{|t|}=\frac{k_{\Bbb
G_n}(0,t\cdot tr(A)e_1+o(t))}{|t|}.$$

Denote by $\hat\kappa_{\Bbb G_n}(0;\cdot)$ the Kobayashi--Buseman
metric of $\Bbb G_n$ at $0,$ that is, the largest norm bounded
above by $\kappa_{\Bbb G_n}(0;\cdot).$ Since $\Bbb G_n$ is a taut
domain, we have that (see \cite{Pang})
$$\lim_{t\to
0}\frac{k_{\Bbb
G_n}(0,t\cdot tr(A)e_1+o(t))}{|t|}=|tr(A)|\hat\kappa_{\Bbb
G_n}(0;e_1).$$
Making use of the equality $\ds\hat\kappa_{\Bbb
G_n}(0;e_1)=\frac{1}{n}$ (cf. \cite{NPZ}), we get that
$$\hat\Omega_n=\{ A \in\mathcal M_n:h_{\hat\Omega_n}(A)=\frac{|tr(A)|}{n}<1\}.$$

\noindent{\bf Remark.} An algebraic approach in the proof of
Proposition \ref{4} also works.
\smallskip

\noindent{\bf Proof of Proposition \ref{nonderent}}. By \eqref{e3},
the equality $\kappa_{\Omega_n}(A;B)$ is equivalent to
$\sigma_{\ast,A}(B)=0$. By property (8) in Proposition \ref{equiv}
and its proof, we have a matrix $Y\in\mathcal M_n$ such that
$-YA+AY=B$. Then the mapping $\lambda\to e^{-\lambda Y} A e^{\lambda
Y}$ satisfies all the required properties.
\smallskip

\noindent{\bf Proof of Conjecture \ref{6} for n=2.} If $A$ is
derogatory, it is scalar. Then $sp(B)=0$ and so the linear mapping
$\lambda\to A+\lambda B$ does the job.

Let $A$ be non-derogatory. Choose $r>0$ and $\varphi\in\mathcal
O(r\Bbb D,\Omega_n)$ such that $\varphi(0)=A,$ $\varphi'(0)=B$ and
$sp(\varphi(\lambda))=sp(A)$ for any $\lambda\in r\Bbb D.$ Then
$$\varphi(\lambda)=A+\lambda B+\lambda^2\psi(\lambda),\quad
\psi\in\mathcal O(r\Bbb D,\mathcal M_n).$$ The Taylor expansion
shows that the condition $sp(\varphi)=sp(A),$ that is
$$tr(\varphi)=tr(A)\hbox{ and }\det\varphi=\det A,$$ is equivalent
to $tr(B)=tr(\psi)=0$ and
$$f(A,B)=\det B+f(A,\psi)=f(B,\psi)=\det\psi=0,$$ where
$$f(C,D)
=c_{11}d_{22}+c_{22}d_{11}-c_{12}d_{21}-c_{21}d_{12},\quad
C,D\in\mathcal M_2$$ Observe that the quadratic mapping
$\lambda\to A+\lambda B+\lambda^2\psi(0)$ satisfies the same
conditions. Therefore it has the desired properties.

\section{Appendix: Proof of Proposition \ref{equiv}}
\label{appendix}

In \cite[definition 3.2.4.1, p. 135]{HoJo}, property (4) is taken
as the definition of {\it nonderogatority}. The fact that (4)
implies (5) is \cite[Theorem 3.2.4.2, p. 135]{HoJo}. The converse
implication is stated in \cite[p. 137]{HoJo}, and proved in
\cite[Corollary 4.4.18, p. 275]{HoJoTop}. The fact that (7) is
equivalent to (4) is part of  \cite[Theorem 4.4.17, p.
275]{HoJoTop}. The equivalence between (1) and (3) is
\cite[Theorem 3.3.15, p. 147]{HoJo}. The equivalence of those
properties with (4) is left as an exercise immediately after this.

From the form of a companion matrix, it is immediate that if we
denote by $e_1$ the first basis vector, its iterates $e_1, Ae_1,
\ldots, A^{n-1}e_1$ generate $\C^n$, and similarity preserves thsi
property. Conversely, if one has a cyclic vector, it is easy to
see that the space is actually generated by the first $n$ iterates
as above, and that they must form a basis, in which the matrix
will take the companion form. So (2) is equivalent to (1).

We now move on to the statements about ranks. First note that
$$
\Phi_A(I+H):=(I+H)^{-1} A (I+H) = A + (-HA+AH) + O(H^2),
$$
so $\dim \Ker((\Phi_A)_{\ast,I_n}) = \dim \mathcal C(A)$, and, by
the rank theorem, (8) is equivalent to (7). The comment about
maximality follows from \cite[Theorem 4.4.17(d), p. 275]{HoJoTop}.

To study $\sigma_{\ast,A}$, first note that for $P\in \mathcal
M_n^{-1}$,
$$
\sigma_{\ast,A} (H) = \sigma_{\ast,P^{-1}AP}(P^{-1}HP),
$$
so $rank(\sigma_{\ast,A})$ is preserved when we pass to a similar
matrix. Thus, if $A$ verifies (1), we may suppose then that it is
a companion matrix. Choose $H=(h_{i,j})$ such that $h_{i,j}=0,$
$1\le j\le n-1$. Then $A+H$ is also a companion matrix, and
$\sigma$ when restricted to that set is a linear map in the last
column. Then the mapping
$$
\sigma_{\ast,A} (H) = (- h_{n,n}, h_{n-1,n}, \ldots, (-1)^{n-1}h_{1,n}),
$$
is onto $\C^n$. So (6) follows by (1).

To complete the proof of Proposition \ref{equiv}, it is enough to
show (6) implies (4).

Given any $\lambda\in\Bbb C$, let $A_\lambda:= A-\lambda I_n$.
Then $P_A(X)= P_{A_\lambda}(X+\lambda)$, so that $\sigma(A) $ is a
polynomial expression (involving the parameter $\lambda$) of the
components of $\sigma(A_\lambda) $.  Therefore $rank
(\sigma_{\ast,A}) \le rank (\sigma_{\ast,A_\lambda})$.

Suppose now that property (6) holds and (4) does not. Let
$\lambda$ be an eigenvalue such that $\dim \Ker (A-\lambda I_n)\ge
2$. Choose a basis of $\C^n$ containing a basis of $ \Ker
(A-\lambda I_n)$.  In this basis, the matrix $ A-\lambda I_n$
transforms into a matrix with at least two columns which are
identically zero, and therefore $\sigma_n(A-\lambda I_n+H)$ is a
polynomial containing only monomials of degree at least $2$ in the
$h_{i,j}.$ This implies that $(\sigma_n)_{\ast,A-\lambda I_n}=0$
and therefore
$$
rank (\sigma_{\ast,A}) \le rank (\sigma_{\ast,A_\lambda I_n}) \le n-1,
$$
which is a contradiction.

\end{document}